\documentclass{article}
\usepackage{amssymb}	% Allows use of AMS's mathematical symbols
\usepackage{amsmath}    % Allows use of AMS's mathematical commands

\usepackage{latexsym}	% Allows use of old latex symbols

\usepackage{epsfig}
\newtheorem{remark}{Remark}[section]

\newcommand{\qed}{\nobreak \ifvmode \relax \else \ifdim\lastskip<1.5em \hskip-\lastskip \hskip1.5em plus0em minus0.5em \fi \nobreak \vrule height0.75em width0.5em depth0.25em\fi}

\def\T{{\rm T}}

\title{An analytic solution to Wahba's problem}
\author {Yaguang Yang\thanks{
NRC, Office of Research, 21 Church Street, Rockville, 20850. Email:
yaguang.yang@verizon.net} and Zhiqiang Zhou\thanks{
NASA Langley Research Center, Hampton, VA 23681. Email: zhiqiang.zhou@nasa.gov}
}
\date{\today}

\begin{document}

\maketitle    % This command generates the title.

\begin{abstract}
All spacecraft attitude estimation methods are based on Wahba's optimization problem. This problem can be reduced to finding the largest eigenvalue and the corresponding eigenvector for Davenport's $K$-matrix. Several iterative algorithms, such as QUEST and FOMA, were proposed, aiming at reducing the computational cost. But their computational time is unpredictable because the iteration number is not fixed and the solution is not accurate in theory. Recently, an analytical solution, ESOQ was suggested. The advantages of analytical solutions are that their computational time is fixed and the solution should be accurate in theory if there is no numerical error. In this paper, we propose a different analytical solution to the Wahba's problem. We use simple and easy to be verified examples to show that this method is numerically more stable than ESOQ, potentially faster than QUEST and FOMA. We also use extensive simulation test to support this claim.
\end{abstract}


%\newpage

\section{ Introduction}
Attitude determination is a crucial part to the success of every spacecraft system \cite{wertz78}.
The spacecraft attitude determination and estimation problem is an optimization problem
formulated by Wahba \cite{wahba65} in 1965. The first solution was provided by Davenport \cite{daven65}.
Davenport proved that Wahba's problem is equivalent to finding the largest eigenvalue and the
corresponding eigenvector for the so-called K-matrix. As the computers and related algorithms at
the time were not powerful enough to use Davenport's method to the real time attitude and control
system of MAGSAT spacecraft \cite{shuster81}, the QUEST algorithm \cite{shuster81} was devised to meet this need.
Since then, QUEST has been widely recognized and used in many spacecraft attitude determination 
and control systems though several other methods, such as ESOQ \cite{mort97} and FOMA \cite{markley93}
were proposed. QUEST and FOMA use Newton's iteration to find the largest 
solution for a quartic polynomial expressed in different ways. Although all flight experiences were 
successful for QUEST method, \cite{markley00} demonstrated that
QUEST may not converge by using a specific example. In fact, it is well known that Newton's method
is inadequate for general use since it may fail to converge to a solution. Even if it does converge,
its behavior may be erratic in regions where the function is not convex \cite{nocedal99}.
Moreover, since the solution of QUEST and FOMA depends on iteration and is numerical which means that
the solution is not precise and the computational time is not predictable.

Noticing that a polynomial of degree 4 (quartics) admits analytic solution, Mortari devised a closed-form 
solution to the Wahba's problem \cite{mort97} where instead of using Newton's method, a set of formulas 
for the roots of quartics is used\footnote{There are several sets of
different analytic formulas for the roots of quartic polynomial \cite{evans94}.}. In this paper, 
we propose an alternative analytic solution which is recently discovered by Shmakov \cite{Shmakov11}. 
The main difference between Mortari's ESOQ and the proposed solution is that the auxiliary cubics  
of the latter is in depressed form while the auxiliary cubics of the former is not. Therefore,
the method based on Shmakov's formulae is more efficient. The simulation test shows that our proposed
method is numerically more stable than ESOQ and comparable to QUEST.

The remaining paper is organized as follows: Section 2 describes Wahba's problem and the solutions
of Davenport, QUEST, ESOQ, and FOMA; followed by our proposed analytic solution. Section 3 provides
some numerical test results to demonstrate the feasibility of the proposed methods. Conclusion 
remarks are summarized in Section 4.

\section{Wahba's problem and solutions}

\subsection{Wahba's problem and numerical solutions}

Given a set of unit vectors $r_i$ that are the
representation of observations of objects in the reference frame, and a set of unit
vectors $b_i$ that are the representation of measurements in the spacecraft body frame,
the attitude of the spacecraft is defined by an orthogonal matrix $A$ that meets the
relation 
\[
Ar_i=b_i.
\]
Therefore, the attitude determination problem is finding the orthogonal matrix
$A$ that minimizes Wahba's loss function
\begin{equation}
L(A)=\frac{1}{2} \sum_{i=1}^{n} a_i | b_i -Ar_i |^2,
\label{wahba}
\end{equation}
where $a_i$, $i=1,\ldots,n$, represents the relative weight associated with the $i$th measurement
error. Solving Wahba's problem is not difficult, but solving it accurately and efficiently
requires some serious effort. All popular methods, such as QUEST \cite{shuster81},
ESOQ \cite{mort97}, and FOMA \cite{markley93}, use Davenport's q-method \cite{daven65,wertz78} 
which converts Wahba's loss function (\ref{wahba}) into a quadratic form
\begin{equation}
\bar{q}^{\T} K \bar{q},
\end{equation}
where $\bar{q}^{\T} = [q^{\T} \cos \left(\frac{\theta}{2} \right), \sin \left(\frac{\theta}{2} \right)]$ 
is a rotational quaternion that brings $r_i$ to $b_i$, where $q$ is the unit vector of rotational axis, $\theta$ is the rotational angle
and $K$ is a $4 \times 4$ matrix given by
\begin{equation}
K= \left[ \begin{array}{cc}
B+B^{\T} - tr[B]I & z \\  z^{\T} & tr[B]
\end{array} \right]
\label{kMatrix}
\end{equation}
where $B=\sum_{i=1}^{n} a_i b_i r_i^{\T}$ is the attitude profile matrix, $I$ is a $3 \times 3$
identity matrix, $z^{\T}=[B_{23}-B_{32}, B_{31}-B_{13}, B_{12}-B_{21}]$, and 
$tr[B]= \sum_{i=1}^{n} B_{ii}$ is the trace of the matrix $B$. By introducing the Lagrange multiplier
$\lambda$ for the unit constraint of $\| \bar{q} \| =1$, Wahba's problem
is reduced to Davenport's problem
\begin{equation}
\max_{\lambda,\bar{q}} \bar{q}^{\T} K \bar{q} - \lambda \bar{q}^{\T} \bar{q}.
\label{obj}
\end{equation}
Taking the derivative of (\ref{obj}) gives the optimal solution which satisfies
\begin{equation}
K \bar{q} = \lambda \bar{q}.
\label{eig}
\end{equation}
Therefore, the optimization problem is reduced to finding the largest eigenvalue of $K$ and its
corresponding eigenvector.

By using the Cayley-Hamilton theorem (cf. \cite{Rugh93}), Shuster \cite{shuster81} derived the first analytic
formula of the characteristic polynomial of the $K$-matrix and suggested using Newton's method to
find the largest eigenvalue $\lambda_{opt}$ iteratively. The optimal attitude $\bar{q}_{opt}$ 
can be obtained by using $\lambda_{opt}$ and the Gibbs vector \cite{shuster81}. This method is
widely recognized and is refereed to as the QUEST method. In 1993, Markley \cite{markley93} derived an 
equivalent characteristic polynomial for the $K$-matrix and also used Newton's method for his 
expression of the polynomial to find the largest eigenvalue $\lambda_{opt}$ iteratively. Markley's method
is now refereed to as the FOMA algorithm.

\subsection{Analytic solution}
Noticing that quartics admits analytical solutions, Mortari \cite{mort97} found a closed-form solution for Wahba's problem. We will show by an example that this method can be numerically unstable.
Our proposed algorithm will use the characteristic polynomial of the $K$-matrix presented in \cite{mort97} which is given as follows.
\begin{equation}
p(x)=x^{4}+ax^3+bx^2+cx+d=0,
\label{char}
\end{equation}
where $a=0$, $b=-2(tr[B])^2+tr[adj(S)]-z^{\T}z$, $S=B+B^{\T}$, $adj(S)$ the adjugate matrix of $S$,
$c=-tr[adj(K)]$, and $d=\det(K)$ are all known parameters. It is well-known that a polynomial of degree $4$
admits analytical solutions. Several different methods were proposed in the last several hundred
years \cite{evans94}. Recently, Shmakov \cite{Shmakov11} found a universal method to find the roots of the 
general quartic polynomial. A special case of this method is simpler than all previous methods and it 
can be directly adopted to solve (\ref{char}). We summarize the steps as follows.

First, (\ref{char}) can be factorized as the product of two quadratic polynomials as
\begin{equation}
(x^2+g_1x+h_1)(x^2+g_2+h_2)=x^4+(g_1+g_2)x^3+(g_1g_2+h_1+h_2)x^2+(g_1h_2+g_2h_1)x+h_1h_2=0.
\label{fact}
\end{equation}
Moreover, $g_1$, $g_2$, $h_1$, and $h_2$ are solutions of two quadratic equations defined by
\begin{subequations}
\begin{align}
g^2-ag+\frac{2}{3}b-y=0 \\
h^2-\left( y+\frac{b}{3} \right) h+d=0
\end{align}
\label{quad}
\end{subequations}
where $y$ is the real root(s) of the following cubic polynomial
\begin{equation}
y^3+py+q=y^3+\left(  ac-\frac{b^2}{3}-4d \right) y +\left(  \frac{abc}{3}
-a^2d-\frac{2}{27}b^3-c^2+\frac{8}{3}bd \right) = 0.
\label{cubic}
\end{equation}
The roots of the cubic equation can be obtained by the famous Cardano's formula \cite{poly07}
\begin{subequations}
\begin{align}
y_1=\sqrt[3]{-\frac{q}{2}+\sqrt{\left( \frac{q}{2} \right)^2 + \left( \frac{p}{3} \right)^3}}
+\sqrt[3]{-\frac{q}{2}-\sqrt{\left( \frac{q}{2} \right)^2 + \left( \frac{p}{3} \right)^3}}
\label{cardanoA} \\
y_2=\omega_1 \sqrt[3]{-\frac{q}{2}+\sqrt{\left( \frac{q}{2} \right)^2 + \left( \frac{p}{3} \right)^3}}
+\omega_2 \sqrt[3]{-\frac{q}{2}-\sqrt{\left( \frac{q}{2} \right)^2 + \left( \frac{p}{3} \right)^3}}
\label{cardanoB} \\
y_3=\omega_2 \sqrt[3]{-\frac{q}{2}+\sqrt{\left( \frac{q}{2} \right)^2 + \left( \frac{p}{3} \right)^3}}
+\omega_1 \sqrt[3]{-\frac{q}{2}-\sqrt{\left( \frac{q}{2} \right)^2 + \left( \frac{p}{3} \right)^3}}
\label{cardanoC}
\end{align}
\label{cardano}
\end{subequations}
where $\omega_1 = \frac{-1+i\sqrt{3}}{2}$ and $\omega_2 = \frac{-1-i\sqrt{3}}{2}$. It is well-known 
that (\ref{cubic}) has either one real solution or three real solutions. If the discriminate
\[
\Delta=\left( \frac{q}{2} \right)^2 + \left( \frac{p}{3} \right)^3 > 0,
\]
then (\ref{cubic}) has a real solution given by (\ref{cardanoA}), and a pair of complex conjugate solutions given
by (\ref{cardanoB}) and (\ref{cardanoC}). If $\Delta = 0$, the (\ref{cubic}) has three zero solutions.
If $\Delta < 0$, then (\ref{cubic}) has three distinct real solutions. In this case, to avoid complex operations, 
the solutions can be given in a different form. Let $r=\sqrt{-\left( \frac{p}{3} \right)^3}$, 
$\theta=\frac{1}{3} \arccos \left( -\frac{q}{2r}  \right)$, then the three real solutions are given by
\begin{subequations}
\begin{align}
y_1= 2 r^{\frac{1}{3}} \cos \left(\theta \right),
\label{1root} \\
y_2= 2 r^{\frac{1}{3}} \cos \left(\theta+\frac{2 \pi}{3} \right),
\label{2root} \\
y_3= 2 r^{\frac{1}{3}} \cos \left(\theta+\frac{4 \pi}{3} \right).
\label{3root}
\end{align}
\end{subequations}
Given a real $y$, 
from (\ref{quad}) or (\ref{1root}), we have 
\begin{subequations}
\begin{align}
g_{1,2}= \pm \sqrt{y-\frac{2}{3}b},
\label{g12} \\
h_{1,2}=\frac{y+\frac{b}{3}\pm \sqrt{(y+b/3)^2-4d}}{2}
\label{h12} \end{align}
\end{subequations}
In view of (\ref{fact}), it is worthwhile to notice that the following relations must be held 
\begin{subequations}
\begin{align}
(g_1+g_2)=a, \label{g12h12a} \\
g_1g_2+h_1+h_2=b, \label{g12h12b} \\
g_1h_2+g_2h_1=c, \label{g12h12c} \\
h_1h_2=d, \label{g12h12d} 
\end{align}
\end{subequations}
where (\ref{g12h12a}), (\ref{g12h12b}), and (\ref{g12h12d}) do not depend on the selections
of $g_1$, $g_2$, $h_1$, and $h_2$ (these relations always hold), but (\ref{g12h12c}) does depend
on the choices of $g_1$, $g_2$, $h_1$, and $h_2$. In practice, we can always assume that $g_1$
takes positive sign in (\ref{g12}) and $g_2$ takes minus sign in (\ref{g12}); we then try that $h_1$
takes positive sign in (\ref{h12}) and $h_2$ takes minus sign in (\ref{h12}); if (\ref{g12h12c})
holds, we get the correct selection; otherwise, $h_1$ takes minus sign in (\ref{h12}) and $h_2$ takes 
positive sign in (\ref{h12}) so that (\ref{g12h12c}) holds.
Finally, the roots of the quartic (\ref{char}) are given by
\begin{subequations}
\begin{align}
x_{1,2}= \frac{-g_1  \pm \sqrt{g_1^2-4h_1}}{2},
\\
x_{3,4}= \frac{-g_2  \pm \sqrt{g_2^2-4h_2}}{2}.
\end{align}
\label{solution}
\end{subequations}

\begin{remark}
It is worthwhile to point out that there is no quadratic term in (\ref{cubic}). Therefore, the
cubics is in depressed form.
\end{remark}

\section{Numerical test}

The proposed analytic method, QUEST-Newton method, ESOQ, and Davenport's q-method have been implemented in MATLAB and tested against each other.

\vspace{0.1in}\noindent
{\bf A simple problem:} The first simple test is the following problem.
\begin{equation}
p(x)=x^{4}+ax^3+bx^2+cx+d=0,
\end{equation}
where $a=0$, $b=-2$, $c=0$, and $d=1$. The problem has two positive solution of $x=1$
and two negative solution of $x=-1$. 
The analytic method finds all solutions without numerical error. Starting from $x=1.1$,
the Newton method finds the largest positive solution $x=1.00000001251746$ after $23$ iterations. This means that QUEST and FOAM may not always converge in one or two iterations as claimed. Their computational time is not fixed and their computation cost may be expensive!

\vspace{0.1in}\noindent
{\bf Test cases in \cite{markley93}:} 
Twelve cases taken from \cite{markley93} are tested to compare the performance of the proposed 
analytic method with the QUEST-Newton, ESOQ, and Davenport's q-method. The observation vectors $b_i$ are calculated from 
the reference vectors $r_i$ with measurement errors $n_i$ by
\[
b_i=Ar_i+n_i.
\]
The true matrix $A$ is given by \cite{markley93}
\[
A=\left[
\begin{array}{ccc}
0.352 & 0.864 & 0.360 \\
-0.864 & 0.152 & 0.480 \\
0.360 & -0.480 & 0.800
\end{array} \right].
\]
The measurement errors are zero-mean Gaussian white noises simulated by 4000 values for each test case. 
The reference vectors and the measurement standard deviations $\sigma_i$ for the twelve test cases 
are listed in Table 1. 

\begin{table}
\begin{center}
\caption{Simulation cases summary}
\begin{tabular}{|r|c|c|}
\hline 

Case & Reference vectors &  Measurement standard deviations \\
\hline
1   & $r_1=[1,0,0]^{\T}, \hspace{0.1in} r_2=[0,1,0]^{\T}, \hspace{0.1in} r_3=[0,0,1]^{\T}$   & $\sigma_1=\sigma_2=\sigma_3=10^{-6}$ rad \\ \hline
2   & $r_1=[1,0,0]^{\T}, \hspace{0.1in} r_2=[0,1,0]^{\T}$   & $\sigma_1=\sigma_2=10^{-6}$ rad  \\ \hline
3   & $r_1=[1,0,0]^{\T}, \hspace{0.1in} r_2=[0,1,0]^{\T}, \hspace{0.1in} r_3=[0,0,1]^{\T}$   & $\sigma_1=\sigma_2=\sigma_3=0.01$ rad \\ \hline
4   & $r_1=[1,0,0]^{\T}, \hspace{0.1in} r_2=[0,1,0]^{\T}$   & $\sigma_1=\sigma_2=0.01$ rad \\ \hline
5   & $r_1=[0.6,0.8,0]^{\T}, \hspace{0.1in} r_2=[0.8,-0.6,0]^{\T}$   & $\sigma_1=10^{-6}$ rad, $\sigma_2=0.01$ rad  \\ \hline
6   & $r_1=[1,0,0]^{\T}, \hspace{0.1in} r_2=[0,0.01,0]^{\T}, \hspace{0.1in} r_3=[0,0,0.01]^{\T}$   & $\sigma_1=\sigma_2=\sigma_3=10^{-6}$ rad  \\ \hline
7   & $r_1=[1,0,0]^{\T}, \hspace{0.1in} r_2=[1,0.01,0]^{\T}$   & $\sigma_1=\sigma_2=10^{-6}$ rad \\ \hline
8   & $r_1=[1,0,0]^{\T}, \hspace{0.1in} r_2=[1,0.01,0]^{\T}, \hspace{0.1in} r_3=[1,0,0.01]^{\T}$   & $\sigma_1=\sigma_2=\sigma_3=0.01$ rad \\ \hline
9   & $r_1=[1,0,0]^{\T}, \hspace{0.1in} r_2=[1,0.01,0]^{\T}$   & $\sigma_1=\sigma_2=0.01$ rad \\ \hline
10  & $r_1=[1,0,0]^{\T}, \hspace{0.1in} r_2=[0.96,0.28,0]^{\T}, \hspace{0.1in} r_3=[0.96,0,0.28]^{\T}$   & $\sigma_1=10^{-6}$ rad, $\sigma_2=\sigma_3=0.01$ rad \\ \hline
11  & $r_1=[1,0,0]^{\T}, \hspace{0.1in} r_2=[0.96,0.28,0]^{\T}$   & $\sigma_1=10^{-6}$ rad, $\sigma_2=0.01$ rad \\ \hline
12  & $r_1=[1,0,0]^{\T}, \hspace{0.1in} r_2=[0.96,0.28,0]^{\T}$   & $\sigma_1=0.01$ rad, $\sigma_2=10^{-6}$ rad \\ \hline
\end{tabular}
\end{center}
\end{table}

The detailed explanation of the test cases is described in \cite{markley93}.  
The simulation results of the estimated errors $\phi$ are listed in Table 2,
where $\phi$ is calculated the same way as suggested in \cite{markley93}
\[
\phi=2\sin^{-1} \left( \frac{ \| A_{e}-A \|}{\sqrt{8}}   \right),
\]
where $A_{e}$ is the estimated rotational matrix. 

\begin{table}
\begin{center}
\caption{ Estimation errors }
\begin{tabular}{|r|c|c|c|c|}
\hline 

Case & Analytic method $\phi$ (deg) &  QUEST method $\phi$ (deg) & ESOQ method $\phi$ (deg)  & q-method  $\phi$ (deg) \\
\hline
1   &  6.495694956077782e-05  &   6.495694956097059e-05 & 1.10609465012e+02 & 6.49569495609e-05 \\ \hline
2   & 8.324164015961696e-05    & 8.324223749065883e-05 & 1.645665577978247 &  8.32422374907e-05 \\ \hline
3   & 0.649531332307863    & 0.649531332307864 & 1.09449329754e+02 & 0.649531332307864 \\ \hline
4   & 0.832408546987256    & 0.832408547860284 & 0.833026446402269  & 0.832408547860284 \\ \hline
5   & 0.557528700788137    & 0.557528701877667 & 0.557646345900295  & 0.557528701877667 \\ \hline
6   & 6.495694956077782e-05   & 6.495694956097059e-05 & 1.10609465012e+02 & 6.49569495609e-05 \\ \hline
7   & 8.324164015961696e-05   & 8.324223749065883e-05 &  1.645665577978247 & 8.32422374907e-05 \\ \hline
8   & 0.649531332307863    & 0.649531332307864 & 1.09449329754e+02 & 0.649531332307864 \\ \hline
9   & 0.832408546987256   & 0.832408547860284 & 0.833026446402269 & 0.832408547860284 \\ \hline
10  & 1.371174492955960   & 1.371174492966333 & 3.628812853257176 & 1.371174492966333 \\ \hline
11  & 1.685838524732360      & 1.685841533993944 & 1.787951534231175 & 1.685841533993952 \\ \hline
12  & 1.670635461315306    & 1.670644941845449 & 1.907061027159859 & 1.670644941845431 \\ \hline
\end{tabular}
\end{center}
\label{testMarkley}
\end{table}
In Table 2, the estimated errors are the means of 4000 runs with the white noise for each 
test case. The results show that for all 12 cases the proposed analytic method is slightly better than the QUEST-Newton and Davenport's q-method. The proposed method is significantly better than ESOQ in cases 1, 3, 6, and 8. Further investigation shows that the largest eigenvalue $\lambda_{opt}$ of the characteristic equation obtained by ESOQ in some case is different from the expected value. One example taken from test case 1 has the following coefficients $a=0$, $b=-0.666666666666667$, $c=-0.296296296294793$, and $d=-0.037037037036536$ in the quartics. The largest eigenvalue given by ESOQ is $\lambda_{opt}=0.333337550713685-i0.000000000001161$, which is not a real number, though it is close to a real number. Moreover, its real part is far away from $\lambda_{opt}=0.999999999999155$ given by other methods.

\section{Conclusions} 
 
We have proposed an analytical solution to Wahba's problem. This solution is accurate
in theory, and the computational time is predictable and fixed. Numerical tests show that
the analytic solution is as stable and accurate as QUEST and q-method but better than ESOQ. Therefore, this solution
may be better than all existing algorithms that have been widely used in aerospace industry. The Matlab code can be obtained from the authors.


\begin{thebibliography}{99}
\bibitem{wertz78} {Wertz, J. (ed.),} {\it Spacecraft Attitude Determination and Control}, 
Kluwer Academic Publishers, Dordrecht, Holland, 1978, pp. 426-428.
\bibitem{wahba65} {Wahba, G.,} ``A least squares estimate of spacecraft attitude,''
{\it SIAM Review}, Vol. 7, No.3, 1965, p.409.
\bibitem{shuster81} {Shuster, M.D., and Oh, S.D.,} ``Three-axis attitude determination from 
vector observations,'' {\it Journal of Guidance and Control}, Vol. 4, No.1, 1981, pp. 70-77.
\bibitem{mort97} {Mortari, D.,} ``{ESOQ}: A closed-form solution to the {W}ahba problem,''
{\it Journal of Astronautical Sciences}, Vol. 45, No. 2, 1997, pp. 195-205.
\bibitem{markley93} {Markley, F.L.,}  ``Attitude determination using vector observations: a fast 
optimal matrix algorithm,'' {\it Journal of Astronautical Sciences}, Vol. 41, No.2, 1993, pp. 261-280.
\bibitem{daven65} {Davenport, P.,} ``A vector approach to the algebra of rotations with applications,''
NASA, X-546-65-437, 1965.
\bibitem{Rugh93} {Rugh, W.J.,} {\it Linear System Theory,} Prentice-Hall, Inc., Englewood Cliffs, 
New Jersey, 1993, pp.4-5.
\bibitem{markley00} {Markley, F.L., and Mortari, D.,} ``Quaternion attitude estimation using 
vector observations,'' {\it Journal of Astronautical Sciences}, Vol. 48, No.2, 2000, pp. 359-380.
\bibitem{nocedal99} Nocedal, J., and Wright, S.J., {\it Numerical Optimization}, Springer-Verlag,
{New York}, 1999, p. 135.
\bibitem{evans94} {Herbison-Evans, D.,} ``Solving Quartics and Cubics for Graphics,'', TR94-487, 
University of Sydney, Sydney, Australia, 1994.
\bibitem{Shmakov11} {Shmakov, S. L.,} ``A universal method of solving quartic equations,''
{\it International Journal of Pure and Applied Mathematics}, Vol. 71, No. 2, 2011, pp. 251-259.
\bibitem{poly07} {Polyanin, A.D., and Manzhirov, A.V.,} {\it Handbook of Mathematics For Engineers 
and Scientists}, Chapman \& Hall/CRC, Noca Raton, FL, 2007, pp. 158-159.
\end{thebibliography}
\end{document}